\documentclass{elsart}
\usepackage{amssymb}
\usepackage{amsfonts}
\usepackage{amsmath}

\input{tcilatex}
\begin{document}

\begin{frontmatter}%

\title{Two deformed Pascal's triangles and its new properties}%

\author[1]{Jishe Feng}
\thanks[1]{Corresponding author. E-mail: gsfjs6567@126.com.}
\author{,Cunqin Shi}
\author{, Huani Zhao}%

\address{Department of Mathematics, Longdong University,  Qingyang,  Gansu,  745000,  China
E-mail: gsfjs6567@126.com.}%

\begin{abstract}%

In this paper, firstly, by a determinant of deformed Pascal's triangle,
namely the normalized Hessenberg matrix determinant, to count Dyck paths, we
give another combinatorial proof of the theorems which are of Catalan
numbers' determinant representations and the recurrence formula. Secondly, a
determinant of normalized Toeplitz-Hessenberg matrix, whose entries are
binomials, arising in power series, we derive new four properties of
Pascal's triangle.

2010 Mathematics Subject Classification: 05A15, 11B65.

\end{abstract}%

\begin{keyword}%
Pascal's triangle, Hessenberg matrix, Catalan number, Dyck path%
\end{keyword}%

\end{frontmatter}%

\section{Introduction}

The so-called Pascal matrix $P=\left( \binom{i}{j}\right) _{i,j\geq 0}$ is
an infinite matrix defined by $p_{ij}=\binom{i}{j}$ for $j\leq i$ and $%
p_{ij}=0$ for $j>i$. Gessel \cite{Gessel} gave a combinatorial
interpretation for the sequential principal minor of the Pascal matrix $P$,
which involves configurations of nonintersecting paths is related to Young
tableaux and hook length formulae. In this paper, we deform Pascal matrix to
a Hessenberg matrix $A$ which is another deformed Pascal's triangle matrix
like this, 
\begin{equation}
A=\left( 
\begin{array}{ccccccc}
\binom{1}{1} &  &  &  &  &  &  \\ 
\binom{2}{0} & \binom{2}{1} & \binom{2}{2} &  &  &  &  \\ 
& \binom{3}{0} & \binom{3}{1} & \binom{3}{2} & \binom{3}{3} &  &  \\ 
&  & \ddots  & \ddots  & \ddots  &  &  \\ 
&  &  & \ddots  & \ddots  & \ddots  &  \\ 
&  &  &  & \binom{n}{0} & \binom{n}{1} & \cdots  \\ 
&  &  &  &  & \cdots  & \cdots 
\end{array}%
\right) .  \label{a1}
\end{equation}%
For convenience, let us denote the $n$th sequential principal minor of order 
$n$ of matrix $A$ by $A_{n}$, namely, $A_{n}=det(a_{ij})_{i,j=0}^{n-1}=\left%
\vert \binom{i+1}{j-i+1}_{i,j=0}^{n-1}\right\vert $.

Catalan number is a ubiquitous sequence of numbers in mathematics. In the
recent book Stanley \cite{Stanley} presents 214 different kinds of objects
that are counted using Catalan numbers. In this paper, we will use the
method of the enumeration of the Dyck paths by a determinant, give a simple
proof of the following two theorems.

\begin{thm}
\cite{vv}For a positive integer $n$, let $C_{n}$ denote the $n$th Catalan
number. Then%
\begin{equation}
C_{n}=A_{n}=\left\vert 
\begin{array}{cccccc}
\binom{1}{1} &  &  &  &  &  \\ 
\binom{2}{0} & \binom{2}{1} & \binom{2}{2} &  &  &  \\ 
& \binom{3}{0} & \binom{3}{1} & \binom{3}{2} & \binom{3}{3} &  \\ 
&  & \ddots & \ddots & \ddots &  \\ 
&  &  & \ddots & \ddots & \ddots \\ 
&  &  &  & \binom{n}{0} & \binom{n}{1}%
\end{array}%
\right\vert .  \label{aa2}
\end{equation}
\end{thm}

\begin{thm}
\cite{Luo} \cite{vv} With the initial condition, $C_{1}=1$, and positive
integer $n>1$. For Catalan number $C_{n}$, there exists a recurrence formula,%
\begin{equation}
C_{n}=\dsum\limits_{j\geq 1}(-1)^{j-1}C_{n-j}\binom{n-j+1}{j}.  \label{a2}
\end{equation}
\end{thm}

This recurrence formula was firstly given by Ming Antu, who was Mongolian
astronomer, mathematician, and topographic scientist \cite{Stanley} \cite%
{Luo}.

Tomislav Do\v{s}li\'{c} \cite{vv} uses counting perfect matchings in a
suitably chosen class of graphs, but as he said that the method comes from a
narrow area. For this reason, we give another combinatorial proof.

Alfred Inselberg \cite{Alfred} gave a Toeplitz-Hessenberg matrix arising in
power series and yielding an asymptotic formula for Bernoulli numbers. In
Section 4, we set the binomials as the entries of the Toeplitz-Hessenberg
matrix and get some new properties of Pascal's triangle.

\section{The Hessenberg determinant and the enumeration method of Dyck paths}

An upper Hessenberg matrix $H_{n}=(h_{ij}),i,j=1,2,\cdots ,n$, is a special
kind of square matrix, such that $h_{i,j}=0$ for $i>j+1$. Ulrich Tamm \cite%
{Tamm} give the concept of the Hessenberg matrix in a normalized form, i.e. $%
h_{i+1,i}=1$ for $i=1,\cdots ,n-1$.\ 
\begin{equation*}
H_{n}=\left( 
\begin{array}{cccccc}
h_{1,1} & h_{1,2} & h_{1,3} & \cdots & \cdots & h_{1,n} \\ 
1 & h_{2,2} & h_{2,3} & \cdots & \cdots & h_{2,n} \\ 
& 1 & h_{3,3} & \cdots & \cdots & h_{3,n} \\ 
&  & \ddots & \ddots &  &  \\ 
&  &  & \ddots & \ddots & h_{n-1,n} \\ 
&  &  &  & 1 & h_{n,n}%
\end{array}%
\right) .
\end{equation*}%
Using the method of Laplace expansion in terms of the last row, we
recursively calculate the determinant of matrix $H_{n}$ as this,%
\begin{equation}
\left\vert H_{n}\right\vert
=\dsum\limits_{i=1}^{n-1}(-1)^{i-1}h_{n-i+1,n}\left\vert H_{n-i}\right\vert .
\label{a3}
\end{equation}

There are a lot of relations between the determinant of Hessenberg matrix
and many well-known number sequences (see\cite{Tamm} \cite{Feng} \cite%
{Yilmaz} and references therein).

Lattice paths are omnipresent in enumerative combinatorics since they can
represent a plethora of different objects. Especially, Dyck paths are
starting at $(0,0)$ and ending to $(n,n)$ with $E=(1,0)$ step and $N=(0,1)$
step that never goes above the line $y=x$ and never goes below the $x$-axis.
It is well-known that the number of Dyck paths equals the nth Catalan number 
\cite{Stanley}.

It is to count paths in a region that is delimited by nonlinear upper and
lower boundaries. Let $a_{1}\leq a_{2}\leq \cdots \leq a_{n}$, and $%
b_{1}\leq b_{2}\leq \cdots \leq b_{n}$ be integers with $a_{i}\geq b_{i}$.
We abbreviate $\mathbf{a}=(a_{1},a_{2},...,a_{n})$ and $\mathbf{b}%
=(b_{1},b_{2},...,b_{n})$. Let $L(0,b_{1})\rightarrow (n,a_{n})$ denote the
set of all lattice paths from $(0,b_{1})$ to $(n,a_{n})$ satisfying the
property that for all $i=1,2,...,n$ the height of the $i$th horizontal step
is in the interval $[b_{i},a_{i}]$. Theorem 10.7.1 in \cite{Krattenthaler}
gives a formula for counting these paths. We restate it as follows.

\begin{thm}
\cite{Krattenthaler} Let $\mathbf{a}=(a_{1},a_{2},...,a_{n})$ and $\mathbf{b}%
=(b_{1},b_{2},...,b_{n})$ be integer sequences with $a_{1}\leq a_{2}\leq
\cdots \leq a_{n}$, $b_{1}\leq b_{2}\leq \cdots \leq b_{n}$, and $a_{i}\geq
b_{i}$, $i=1,2,...,n$. The number of all paths from $(0,b_{1})$ to $%
(n,a_{n}) $ satisfying the property that for all $i=1,2,...,n$ the height of
the $i$th horizontal step is between $b_{i}$ and $a_{i}$ is given by%
\begin{equation*}
|L((0,b_{1})\rightarrow (n,a_{n}):\mathbf{b}\leq y\leq \mathbf{a})|=\underset%
{1\leq i,j\leq n}{det}\left( \binom{a_{i}-b_{j}+1}{j-i+1}\right) .
\end{equation*}
\end{thm}

\section{Proofs of Theorem 1 and 2}

\begin{pf}
Set $b=(0,0,\cdots ,0)$ and $a=(0,1,2,3,\cdots ,n)$. Applying Theorem 3, we
get the number of Dyck paths which start at $(0,0)$ to $(n,n)$ satisfying
the property that for all $i=0,1,2,\cdots ,n$ the height of the $i$th
horizontal step is between $0$ and $i$ is given by the determinant (\ref{aa2}%
), $A_{n}$, which equals the $n$th Catalan number $C_{n}$.

Assume the determinant $|A_{0}|=1$, and apply the formula (\ref{a3}) to (\ref%
{aa2}), we can get%
\begin{eqnarray*}
A_{n} &=&\binom{n}{1}A_{n-1}+(-1)\binom{n-1}{2}\binom{n}{0}A_{n-2}+(-1)^{2}%
\binom{n-2}{3}\binom{n-1}{0}\binom{n}{0}A_{n-2} \\
&&+(-1)^{3}\binom{n-3}{4}\binom{n-2}{0}\binom{n-1}{0}\binom{n}{0}%
A_{n-3}+\cdots \\
&=&\dsum\limits_{j\geq 1}(-1)^{j-1}A_{n-j}\binom{n-j+1}{j}.
\end{eqnarray*}%
Because $|A_{0}|=C_{0}=1$, one obtains the formula (\ref{a2}).
\end{pf}

\section{The determinant of Toeplitz-Hessenberg matrix arising in power
series}

In this section, we give a determinant of the Toeplitz-Hessenberg matrix
arising in power series. The new four properties of the determinant in which
entries are binomials are derived. Assume $n$ is a positive integer, as in 
\cite{Alfred}, we set $a_{i}=\binom{n}{i}$ in the determinant $J_{n}$ get a
determinant as%
\begin{equation}
J_{n,m}=\left\vert 
\begin{array}{cccccc}
\binom{n}{1} & 1 &  &  &  &  \\ 
\binom{n}{2} & \binom{n}{1} & 1 &  &  &  \\ 
\vdots & \ddots & \ddots & \ddots &  &  \\ 
\binom{n}{n} &  & \ddots & \ddots & \ddots &  \\ 
& \ddots &  & \ddots & \ddots & 1 \\ 
&  & \binom{n}{n} & \cdots & \binom{n}{2} & \binom{n}{1}%
\end{array}%
\right\vert _{m\times m},  \label{a4}
\end{equation}%
where $m\geq n$ and assume $J_{0,0}=1$, which arises in the power serises of
the reciprocal of the function%
\begin{equation*}
f(x)=(1+z)^{n}=\dsum\limits_{i=0}^{n}\binom{n}{i}z^{i}.
\end{equation*}

One obtains%
\begin{equation*}
g(x)=\frac{1}{f(x)}=\dsum\limits_{j=0}^{\infty
}(-1)^{j}M_{n,j}z^{j}=\dsum\limits_{j=0}^{\infty }(-1)^{j}\binom{n+j-1}{j}%
z^{j},
\end{equation*}

and the following four properties

(a) For $k=0,1,2,\cdots $, the $k$ order sequential principal minor $M_{n,k}$
of $J_{n,m}$ equals 
\begin{equation}
M_{n,k}=\binom{n+k-1}{k}.  \label{a6}
\end{equation}

(b) The series $\{M_{n,k}\}_{0}^{\infty }$ is on the nth column in the
Pascal's triangle. For example, $n=8$, assume that $M_{0}=1,$ $%
\{M_{8,k}\}_{0}^{\infty }\ $is the $8$-th column or row.%
\begin{equation}
M_{8,k}=\left\vert 
\begin{array}{cccccc}
\binom{8}{1} & 1 &  &  &  &  \\ 
\binom{8}{2} & \binom{8}{1} & 1 &  &  &  \\ 
\vdots & \ddots & \ddots & \ddots &  &  \\ 
\binom{8}{8} &  & \ddots & \ddots & \ddots &  \\ 
& \ddots &  & \ddots & \ddots & 1 \\ 
&  & \binom{8}{8} & \cdots & \binom{8}{2} & \binom{8}{1}%
\end{array}%
\right\vert  \label{a9}
\end{equation}

\begin{tabular}{|l|l|l|l|l|l|l|l|l|l|l|l|l|l|l|l|l|}
\hline
1 & 1 & 1 & 1 & 1 & 1 & 1 & \textbf{1} & 1 & 1 & 1 & 1 & 1 & 1 & 1 & 1 & $%
\cdots $ \\ \hline
1 & 2 & 3 & 4 & 5 & 6 & 7 & \textbf{8} & 9 & 10 & 11 & 12 & 13 & 14 & 15 & 
&  \\ \hline
1 & 3 & 6 & 10 & 15 & 21 & 28 & \textbf{36} & 45 & 55 & 66 & 78 & 91 & 105 & 
&  &  \\ \hline
1 & 4 & 10 & 20 & 35 & 56 & 84 & \textbf{120} & 165 & 220 & 286 & 364 & 455
&  &  &  &  \\ \hline
1 & 5 & 15 & 35 & 70 & 126 & 210 & \textbf{330} & 495 & 715 & 1001 & 1365 & 
&  &  &  &  \\ \hline
1 & 6 & 21 & 56 & 126 & 252 & 462 & \textbf{792} & 1287 & 2002 & 3003 &  & 
&  &  &  &  \\ \hline
1 & 7 & 28 & 84 & 210 & 462 & 924 & \textbf{1716} & 3003 & 5005 &  &  &  & 
&  &  &  \\ \hline
\textbf{1} & \textbf{8} & \textbf{36} & \textbf{120} & \textbf{330} & 
\textbf{792} & \textbf{1716} & \textbf{3432} & \textbf{6435} &  &  &  &  & 
&  &  &  \\ \hline
1 & 9 & 45 & 165 & 495 & 1287 & 3003 & \textbf{6435} &  &  &  &  &  &  &  & 
&  \\ \hline
1 & 10 & 55 & 220 & 715 & 2002 & 5005 & \textbf{11440} &  &  &  &  &  &  & 
&  &  \\ \hline
1 & 11 & 66 & 286 & 1001 & 3003 & 8008 & \textbf{19448} &  &  &  &  &  &  & 
&  &  \\ \hline
1 & 12 & 78 & 364 & 1365 & 4368 & 12376 & \textbf{31824} &  &  &  &  &  &  & 
&  &  \\ \hline
1 & 13 & 91 & 455 & 1820 &  &  &  &  &  &  &  &  &  &  &  &  \\ \hline
1 & 14 & 105 & 560 &  &  &  &  &  &  &  &  &  &  &  &  &  \\ \hline
1 & 15 & 120 &  &  &  &  &  &  &  &  &  &  &  &  &  &  \\ \hline
1 & 16 &  &  &  &  &  &  &  &  &  &  &  &  &  &  &  \\ \hline
1 &  &  &  &  &  &  &  &  &  &  &  &  &  &  &  &  \\ \hline
$\vdots $ &  &  &  &  &  &  &  &  &  &  &  &  &  &  &  &  \\ \hline
\end{tabular}

(c) For $k=0,1,2,\cdots $, there are the following relations 
\begin{equation}
\dsum\limits_{i=0}^{k}(-1)^{i}\binom{n+i-1}{i}\binom{n}{k-i}=0.  \label{a7}
\end{equation}

(d) There are two relations, 
\begin{equation}
J_{n,m}=\dsum\limits_{h=1}^{n}(-1)^{h+1}\binom{n}{h}J_{n,m-h}=\binom{n+m-1}{m%
}.  \label{a8}
\end{equation}%
\begin{equation}
\dsum\limits_{h=1}^{n}(-1)^{h+1}\binom{n}{h}\binom{n+m-h-1}{m-h}=\binom{n+m-1%
}{m}  \label{a88}
\end{equation}


\begin{thebibliography}{9}
\bibitem{Alfred} Alfred Inselberg. On determinants of Toeplitz-Hessenberg
matrices arising in power series. Journal of Mathematical Analysis and
Applications 63(1978): 347-353.

\bibitem{Feng} Jishe Feng. Hessenberg Matrices On Fibonacci And Tribonacci
Numbers. Ars Combinatoria, 127(2016) 117-124.

\bibitem{Gessel} I. Gessel, Binomial determinants, paths, and hook length
formulae, Advances in Mathematics 58(1985)300-321.

\bibitem{Krattenthaler} C. Krattenthaler. lattice path enumeration. in
Handbook of Enumerative Combinatorics (M. Bona, ed.), ch.10, pp. 589-678,
Boca Ranton: CRC Press, Boca Raton-London-New York, 2015.

\bibitem{Luo} J. Luo, Ming Antu and his power series expansions. in Sek.
founder of modern mathematics in Japan Springer, Tokyo, 2013, 299-310.

\bibitem{Stanley} R.P. Stanley, Catalan numbers, Cambridge University Press,
Cambridge, 2015.

\bibitem{Tamm} Ulrich Tamm. The determinant of a Hessenberg matrix and some
applications in discrete mathematics. https:// www.math.uni-bielefeld.ed/
ahlswede/pub/tamm/hessen.ps.

\bibitem{vv} Tomislav Do\v{s}li\'{c}. Perfect Matchings, Catalan Numbers,
and Pascal's Triangle, Mathematics Magazine, 80:3 (2007) 219-226.

\bibitem{Yilmaz} F. Yilmaz, D. Bozkurt. Hessenberg matrices and the Pell and
Perrin numbers. Journal of Number Theory, 131 (2011) 1390--1396.
\end{thebibliography}
\end{document}